\input amstex.tex

\documentstyle{amsppt}
\magnification \magstephalf
\pagewidth{6.2in}
\pageheight{8.0in}
\NoRunningHeads
\def\makefootline{\baselineskip=49pt \line{\the\footline}}

\define\la{\lambda}
\define\al{\alpha}

\define\ep{\varepsilon}

\
\define\sign{\operatorname{sign}}
\define\trace{\operatorname{trace}}
\define\leqK{\ \leq\Sb \\ \\ {\!\!\!\!\! K}\endSb\ }
\define\leqk{\ \leq\Sb \\ \\ {\!\!\!\!\! k}\endSb\ }
\define\leqkK{\ \leq\Sb \\ \\ {\!\!\!\!\! k,\,K}\endSb\ }

\topmatter
\title
A negative mass theorem for  the $2$-Torus
\endtitle
\author  K. Okikiolu
\endauthor
\thanks{Date: Revised July 28, 2008.   MSC classes: 58J50, 34A26, 35J60,
53C20.\newline The author would like to thank the University of
Pennsylvania for their hospitality.
 \newline The author was supported  by the National Science Foundation
\#DMS-0302647. } \endthanks
\abstract Let $M$ be a closed surface. For a metric $g$ on $M$,
denote the area element by $dA$ and the Laplace-Beltrami operator by
$\Delta=\Delta_g$.  We define the Robin mass $m(p)$ at the point
$p\in M$ to be the value of the Green function $G(p,q)$ at $q=p$
after the logarithmic singularity has been subtracted off, and we
define $\trace\Delta^{-1}=\int_M m(p)\,dA$. This regularized trace
can also be obtained by  regularization of the spectral zeta
function and is hence a spectral invariant. Furthermore,
$(\trace\Delta^{-1})/A$ is a non-trivial analog for closed surfaces
of the ADM mass for higher dimensional asymptotically flat
manifolds.  We define the $\Delta$-mass of $(M,g)$ to equal
$(\trace\Delta_g^{-1}-\trace \Delta_{S^2,A}^{-1})/A$, where
$\Delta_{S^2,A}$ is the Laplacian on the round sphere of area $A$.
 In this paper we show that in each conformal class $\Cal C$ for the $2$-torus,
there exists a  metric  with negative $\Delta$-mass.  From this it
follows that the minimum of  the $\Delta$-mass on $\Cal C$ is
negative and attained by some metric $g\in\Cal C$.  For this
minimizing metric $g$, one gets a sharp logarithmic
Hardy-Littlewood-Sobolev inequality and an Onofri-type inequality.
We remark that if the flat metric in $\Cal C$ is sufficiently long
and thin then the minimizing metric $g$ is non-flat.
 The proof of our result depends on analyzing the ordinary
differential equation $\phi''=1-e^\phi$ which is equivalent to
$h''=1-1/h$. The solutions are periodic and we  need to establish
quite delicate, asymptotically sharp inequalities relating the
period to the maximum value.
\newline

\noindent
\endabstract
\endtopmatter

\noindent{\bf Section 1. Introduction, Main Results and Summary of
the Proof}.
\medskip

Let $M$ be a smooth, closed, compact surface with a (Riemannian)
metric $g$. Denote the area element of $g$ by $dA$ and the area by
$A$. Let $\Delta=\Delta_g$ denote the Laplace-Beltrami operator for
$g$, given in local coordinates $(x_1,\dots,x_n)$ by
$$
\Delta \ =\ -\sum_{i,j}\ \frac1{\sqrt{\det g}}\
\frac{\partial}{\partial{x_i}} \ \sqrt{\det g}\ g^{ij}
\frac{\partial }{\partial{x_j}} .\tag1.1
$$
The  kernel of $\Delta$ is the constants. Let $\Delta^{-1}$ denote
the inverse operator
$$
\Delta^{-1}\Delta f\ =\ f\ -\ \frac1A\int_M f\,dA.
$$
The Green function $G(p,q)$ for $\Delta$ is the smooth function on
$M\times M\setminus \{(p,p):p\in M\}$ which satisfies
$$
\Delta^{-1} f(p)\ =\ \int_M G(p,q) f(q)\,dA(q).
$$
Denoting the distance from $p$ to $q$ in the metric $g$ by $d(p,q)$, the function
$G(p,q)$ is smooth away from the diagonal  and has an expansion at the diagonal of
the form
$$
G(p,q)\ =\ -\frac1{2\pi}\log d(p,q)\ +\ m(p)\ +\ o(d(p,q)).\tag 1.2
$$
We call the value $m(p)=m_g(p)$ the {\it Robin mass at the point
$p$}.  For a smooth function $\phi$ on $M$, write $A_\phi$ for the
area of $M$ in the metric $e^\phi g$, so

\centerline{$\dsize{ A_\phi\ =\ \int_M e^\phi\,dA.}$}

\proclaim{Conformal change of the Robin mass} If $\phi$ is a smooth
function on $M$ then
$$
m_{e^\phi g}(p)\ =\  m _g(p)\ +\ \frac{\phi}{4\pi}\ -\  \frac{2}{ A_\phi}
(\Delta^{-1}_g e^\phi)(p) \ +\  \frac1{A_\phi^2}\int_M e^\phi
\Delta_g^{-1}e^\phi\,dA. \tag 1.3
$$
\endproclaim
\noindent For the proof, see for example [S1], [S2], [M2] or [O2].
We define
$$
\trace \Delta^{-1}_g \ =\ \int_M m_g(p)\,dA(p).
$$
 This is  a spectral invariant for $\Delta$, since
 it can be
 obtained from the spectral zeta function associated to $\Delta$, see [S1],
 [S2], [M3], or [O2].

 \remark{Remark}
Writing $K(p)$ for the Gaussian curvature of $g$ at $p$, it is shown
in [S1], [S2],  that for any metric $g$ on the $2$-sphere, we have
$$
m_g(p)\ -\ \frac1{2\pi}\Delta^{-1}K(p)\ =\ \frac1A
\trace\Delta^{-1}_g.\tag 1.4
$$
 The left hand side (and hence the right hand side)
 is a $2$-sphere analog  of the
 ADM mass from general relativity.  Indeed, the (Riemannian) ADM mass is
  defined for  asymptotically
 flat manifolds. However, if $M$ is a  compact Riemannian manifold
 of dimension greater than $2$, with positive
 conformal Laplacian, then given a point $p\in M$ we can define
 a mass at $p$ by blowing up the  metric around $p$ using the Green function
 for
 the conformal Laplacian, and taking the ADM mass of the resulting  asymptotically flat metric.
 This amounts to taking the constant term in the asymptotic expansion of the Green function
 for the conformal Laplacian around the point $p$.
 The left hand side of (1.4) is
 the natural non-trivial analog of this for the $2$-sphere.
  Formula (1.4) does not hold for surfaces of higher
 genus.  The left hand side is no longer pointwise constant
 and its fluctuation does not
 have obvious geometric significance.  Therefore we
  consider the right hand side of (1.4) as a natural
  non-trivial
analog of the
 ADM mass for compact surfaces.
\endremark

Now (1.3) immediately gives the following formula, see also [M1].

 \proclaim{Conformal change of
$\trace\Delta^{-1}$ (Morpurgo's Formula)}  If $\phi$ is a smooth
function on $M$, then
$$
\trace \Delta^{-1}_{e^\phi g}\ =\ \int_M m_g e^\phi \,dA\ +\ \frac1{4\pi}\int_M
\phi\, e^\phi \,dA\ -\ \frac1{A_{\phi}} \int_M e^\phi \Delta^{-1}_g e^\phi\,dA. \tag 1.5
$$
\endproclaim
On the round sphere, the right hand side of (1.5) occurs in the
logarithmic Hardy-Littlewood-Sobolev inequality.

 \proclaim{\bf Sharp logarithmic Hardy-Littlewood-Sobolev
inequality on the $S^2$} If $g$ is a round metric on $S^2$ of area
$A$,
$$
 \frac1{4\pi}\int_{S^2} \phi\, e^\phi \,dA\ -\ \frac1A\int_{S^2}\, e^\phi  \Delta^{-1} e^\phi \,dA \ \geq\ 0
$$
holds for all  functions $\phi :S^2\to \Bbb R$ with $ \int_{S^2}
e^\phi\,dA=A$ such that $\int_{S^2} \phi\,e^\phi\,dA$ is finite.
Moreover equality is attained exactly when $e^\phi$ is the Jacobian
of a conformal transformation of $S^2$.
\endproclaim
\medskip
For the proof, see  [On], [CL], [B]. Combining this with (2),
Morpurgo obtained the following. \proclaim{Spectral interpretation
of the logarithmic HLS inequality} Among all metrics on the
$2$-sphere of area $A$, the round metric attains the minimum value
of $\trace\Delta^{-1}$.
\endproclaim

\noindent The behavior of $\trace\Delta^{-1}$ for non-flat metrics
on the torus was first considered in [M1]. Suppose $g_0$ is any flat
metric of unit area on the $2$-torus, and let $\la_1(g_0)$ denote
the lowest eigenvalue of the Laplace-Beltrami operator for $g_0$.
Let $\Cal C_1$ denote the class of metrics conformal to $g_0$ having
unit area.  It was shown in [M1] that if $\la_1(g_0)>8\pi$, then
$g_0$ is a {\it local} minimum for $\trace\Delta^{-1}$ on $\Cal
C_1$.  In [LL1], [LL2], this was improved to a global result in most
cases. Indeed,  it was shown that $g_0$ minimizes
$\trace\Delta^{-1}$ on $\Cal C_1$ provided $\la_1(g_0)\geq \pi^3$,
or $g_0$ is rectangular and $\la_1\geq 8\pi$.   It is well
understood that $g_0$ cannot minimize $\trace\Delta^{-1}$ on $\Cal
C_1$ when  $\la_1(g_0)$ is small. Indeed, it can be observed from
the Kronecker limit formula that when $\la_1(g_0)$ is small, the
value of $\trace\Delta^{-1}$ for $g_0$ is greater than the value for
the round sphere of unit area, as was pointed out in  [DS2].
However, by blowing a spherical bubble, one can construct a family
of  metrics in $\Cal C_1$   for which $\trace\Delta^{-1}$ approaches
the value for the round sphere  (see [O2], [DS2] for different
approaches to this). In this paper, we show that if $T$ is a flat
torus of unit area with $\la_1(T)<8\pi$, then the minimum value of
$\trace\Delta^{-1}$ among conformal metrics of unit area is attained
by a non-flat metric. Although we do not identify this minimizing
metric explicitly, we do construct a candidate, which is
approximately spherical except for a short wormhole joining the
poles.

\proclaim{Theorem 1} Let $T$ be a $2$-dimensional torus with metric
$g_0$. Then  there exists a  metric $g$ in the same conformal class
as $g_0$ and having the same area $A$, such that the Robin mass
$m(x)$ for  $g$ is constant, and strictly less than the Robin mass
for the round sphere of area $A$.
\endproclaim

 This leads to the following result.

 \proclaim{Theorem 2}  Let $T$ be a $2$-dimensional torus with metric $g_0$.  Then among
  metrics in the same conformal class as $g_0$ and having the same
 area $A$, there exists a metric $g$
 which attains the minimum value of
 $\trace\Delta^{-1}$.  Moreover $g$ has constant Robin mass $m(x)$,
 and this is less than the Robin mass of the  round sphere of area $A$.
\endproclaim
\medskip

We remark that if $g_0$ is flat with $\la_1(g_0)<8\pi$, then the
metric $g$ is not flat and the Robin mass for $g$ is less than that
for $g_0$.
\medskip

\proclaim{Corollary 3} (Analogs of Logarithmic HLS inequality and
Onofri's Inequality for the torus.)  For the minimizing metric $g$
of Theorem 2, we have
$$
 \frac1{4\pi}\int_{T} \phi\,e^\phi \,dA\ -\ \frac1A\int_{T}\, e^\phi \Delta^{-1} e^\phi \,dA \ \geq\ 0
\tag 1.6
$$
for all  functions $\phi:T\to\Bbb R$ with $ \int_{T} e^\phi \,dA=A$
such that $\int_{T} \phi\,e^\phi \,dA$ is finite.  Here, $dA$ and
$\Delta$ are associated to $g$. Moreover, for $\phi\in C^\infty(T)$,
$$
\frac{1}{16\pi}\int_T \phi \Delta \phi\,dA\ -\ \log\left(
\frac1A\int_T e^{\phi}\,dA\right)\ +\ \frac1A\int_T \phi\,dA\ \geq\
0.
$$
\endproclaim

To deduce Theorem 2 from Theorem 1, we appeal to Theorem 1 of [O2],
which states that the minimum value of $\trace\Delta^{-1}$ among
metrics conformal to $g_0$ having the same area is attained,
provided there exists a metric conformal to $g_0$ for which the
value of $\trace\Delta^{-1}$ is lower than the value for the round
sphere of the same area.  The proof of that result is a variational
argument very similar in spirit to the proof of the Yamabe theorem
in the non-positive case.  One is trying to find $\phi$ to minimize
(1.5). First one modifies the equation to break the lack of
compactness by replacing $\Delta^{-1}$ in the integral on the right
by $\Delta^{-1-\ep}$.  One can construct a minimizer for the
resulting functional, and one wants this minimizer to converge to a
limit
 as $\ep\to 0$. It is here that one uses the fact that the
value of $\trace\Delta^{-1}$ is lower than that for the round
sphere, which is what prevents bubbles from forming and ensures the
existence of a convergent subsequence as $\ep\to 0$. To deduce
Corollary 3 from Theorem 2, we appeal to Theorem 3 in [O2], which is
just an explicit formulation of the duality between the logarithmic
Sobolev inequality and the Onofri inequality.  For some related
results, see [Ch],  [M2], [M3], [O1], [OsPS1], [S2]. For a
probabilistic interpretation of $\trace\Delta^{-1}$, see [DS1].
\medskip\medskip\medskip

\noindent{\bf Proof of Theorem 1.}  We will quickly show that our
result is related to the  problem of establishing somewhat delicate
inequalities  between the period and the maximum value of solutions
to the ordinary differential equation $\phi''=1-e^\phi$. These
inequalities are established by making just the right Taylor
expansion of the integral formula for the period.

We first remark that under scaling by a constant $e^\la$, the Robin
mass scales as
$$
m_{e^\la g}(p)\ =\  m _g(p)\ +\ \frac{ \la}{4\pi}.
$$
Hence if we can prove the Theorem for area $A=1$, it follows for
arbitrary values of $A$. Furthermore, by the classical
Uniformization Theorem we can assume that $g_0$ is a flat metric on
$T$ with area $1$, and we seek the metric $g=e^\phi g_0$ of area
$1$.
 From (1.3), the condition that the mass $m_{e^\phi g_0}(p)$ is constant is
$$
 \phi\ -\  8\pi (\Delta_0^{-1} e^\phi)\ \text{ is constant},
$$
where $\Delta_0$ is the Laplacian for $g_0$. Applying $\Delta_0$  we
find that this is equivalent to
$$
\Delta_0 \phi\ =\ 8\pi (e^\phi-1). \tag 1.7
$$
We remark  that if $\phi$ satisfies this condition then the metric
$e^\phi g_0$ automatically has area $1$, since
$$
0\ =\ \int_T \Delta_0 \phi\,dA\ =\ 8\pi \int_T (e^\phi-1)\,dA_0.
$$
where $dA_0$ is the area element for $g_0$. We assume that $\phi$
satisfies (1.7). Then (1.5) gives
$$
\trace \Delta^{-1}_{e^\phi g_0}\ =\ \trace\Delta^{-1}_0\ +\
\frac1{8\pi}\int_T \phi (1+e^\phi)\,dA_0, \tag 1.8
$$

 Now we work on a torus with flat metric $g$ of
area $1$, given by  $\Bbb C/\Lambda$ where $\Lambda$ is the lattice
generated by $1/b$ and $a+ib$.  A fundamental domain for the torus
is given by
$$
\left\{x+iy :0\leq y\leq b,\qquad \frac{ay}b\leq x\leq
\frac{ay+1}b\right\}. \tag 1.9
$$
It is a fact that every metric on the torus is conformal to such a
flat metric, with
$$
b\ \geq\ \left(\frac34\right)^{1/4}\ =\ 0.9306... \tag 1.10
$$
For the flat metric $g$ on this torus, we compute in the appendix
using the first Kronecker limit formula that setting
$$
\beta\ =\ \sqrt{\pi}\,b,\tag 1.11
$$
we have
$$
 \trace\Delta_{0}^{-1}\ -\ \trace\Delta_{S^2,1}^{-1}\ =\
\frac1{4\pi} \left( \frac{\beta^2}3\ -\ \log (4\beta^2)\ +\ 1\  -\
4\sum_{n=1}^\infty \log\left| 1-e^{-2n(\beta^2-i\sqrt{\pi}\,\beta
a)}\right| \right),\tag 1.12
$$
where $\Delta^{-1}_{S^2,1}$ is the Laplacian on the round $2$-sphere
of area $1$, see also [Chiu], [S1], [S2]. From this we see that
$$
 \trace\Delta_{0}^{-1}\ -\ \trace\Delta_{S^2,1}^{-1}\ \leq\
\frac1{4\pi} \left( \frac{\beta^2}3\ -\ \log (4\beta^2)\ +\ 1\  -\
4\sum_{n=1}^\infty \log\left| 1-e^{-2n\beta^2}\right|\right). \tag
1.13
$$
From this point, the proof involves some simple numerical
evaluations as well as exact formulas and asymptotic estimates. It
is a fact first pointed out in [DS2] that  that the left hand side
of (1.13) is negative when $\beta$ is small.  To see this, note that
$$
- 4\sum_{n=1}^\infty \log\left| 1-e^{-2n\beta^2}\right|
$$
is decreasing in $\beta$ and is thus bounded by the value at the
endpoint $\beta=\pi^{1/2}(3/4)^{1/4}$, which is
$$
 -4\sum_{n=1}^\infty
\log\left| 1-e^{- 3^{1/2}\pi n }\right|\ <\ 0.02.
$$
On the other hand,
$$
\frac{\beta^2}3\ -\ \log (4\beta^2)\ +\ 1
$$
is convex on the interval $[\pi^{1/2}(3/4)^{1/4}, 2.6]$, and hence
is bounded above there by $-0.04$. Adding these terms, we find that
the right hand side of (1.13) is negative when $\beta\leq 2.6$.
 We see
then that in this case the flat metric $g=g_0$ satisfies the
conclusion of Theorem 1. We only need prove Theorem 1 when
$\beta>2.6$. Noting that $2.6>\pi/\sqrt{2}$, we now complete the
proof of Theorem $1$, by explaining how to find $g$ in the case
$\beta>\pi/\sqrt{2}$.
\medskip

\noindent{\bf Remark}. If $b>1$, then the length of the shortest
geodesic is $1/b$ and the lowest eigenvalue of the Laplace-Beltrami
operator is $\la_1=4\pi^2/b^2=4\pi^3/\beta^2$, so the value
$\beta=\pi/\sqrt{2}$ corresponds to $\la_1=8\pi$.  The value
$\beta=2$ corresponds to $\la_1=\pi^3$.  We remark that when
$\beta\leq 2$, it is shown  in [LL1]  that the flat metric minimizes
$\trace\Delta^{-1}$.  Since the minimum must beat the round sphere,
this again confirms for the case $\beta\leq 2$ , that (1.13) is
negative.
\medskip

Assuming $\phi$ satisfies (1.7), combining (1.8) and (1.13) gives
$$
\multline
\trace \Delta^{-1}_{e^\phi g_0}\ -\ \trace\Delta_{S^2,1}^{-1}\\ \leq\
\frac1{4\pi} \left(\frac1{2}\int_T \phi(1+e^\phi)\,dA_0\ +\  \frac{\beta^2}3\ -\ \log (4\beta^2)\ +\ 1\
-\ 4\sum_{n=1}^\infty \log\left| 1-e^{-2n\beta^2}\right| \right).
\endmultline
\tag 1.14
$$
 We will find $\phi\in C^\infty(T)$ satisfying (1.7) such that $\phi(x+iy)$
 is a function of  $y$ alone, and the right hand side of (1.14) is negative.
We can recast (1.7) and (1.14) in terms of the single variable $y$
so that  Theorem 1 follows from the following:

\proclaim{Theorem 1$'$} For each $b> (\pi/2)^{1/2}$, there exists a
smooth function  $\phi\in C^\infty(\Bbb R)$ satisfying
$$
\align & \frac{d^2\phi}{dy^2}\  =\ 8\pi (1-e^\phi), \tag 1.15\\
& \phi(y+b)=\phi(y)\ \text{ for every }y\in\Bbb R,\tag 1.16\\
& \phi \text{ attains its maximum value }\phi_0\text{ at } y=0,\tag 1.17
\endalign
$$
and such that writing $\beta=\pi^{1/2}b$, we have
$$
\frac1{2b}\int_0^b \phi(1+e^\phi)\,dy\ +\  \frac{\beta^2}3\ -\ \log (4\beta^2)\ +\ 1\
-\ 4\sum_{n=1}^\infty \log\left| 1-e^{-2n\beta^2}\right|\ <\ 0. \tag 1.18
$$
\endproclaim

\remark{Remarks} 1. The condition (1.17) is just thrown in to
eliminate the degree of freedom given by translation invariance. In
fact we choose $\phi$ to have smallest period $b$, which together
with (1.15) and (1.17) determines $\phi$ uniquely.

2. In proving Theorem 1$'$, we will establish a relationship between
the maximum value $\phi_0$ of $\phi$  and the period $b$. A
simplified version is that there exist $\ep_1,\ep_2>0$ such that
$$
e^{\phi_0} + \log 4 + \ep_1 e^{-\phi_0}\ \leq\ \pi b^2\ \leq\
e^{\phi_0} + \log 4 + \ep_2 e^{-\phi_0},\qquad\qquad\text{ for
}b\geq \left(\frac\pi2\right)^{1/2}.
$$
The precise version is that there exist $\ep_1,\ep_2>0$ such that
$$
e^{\phi_0}-\phi_0 + \log 4 + \ep_1 e^{-\phi_0}\ \leq\ \pi b^2 -
\log(\pi b^2)\ \leq\ e^{\phi_0}-\phi_0 + \log 4 + \ep_2
e^{-\phi_0},\qquad\qquad\text{ for }b\geq
\left(\frac\pi2\right)^{1/2}.\tag 1.19
$$

3. In [DS2], conformal factors were chosen for  long skinny flat
tori of area $1$, so that as the length of the flat torus tends to
infinity, the Robin mass of the new metric converges to that of the
round sphere.  From [O2], one sees this can easily be accomplished
by conformal factors which concentrate at a point, but the conformal
factors in [DS2] depend only on the length variable $y$.  In this
paper  we  choose conformal factors which minimize the Robin mass
among one-variable candidates, yielding optimal metrics which beat
the mass of the sphere on every torus. It is unknown whether our
conformal factors give the true minimizer in any case.
\endremark
\medskip

 The rest of the paper is dedicated to proving Theorem 1$'$.
We begin by giving a summary of the proof, and then supply the
details,
\medskip

\noindent{\bf Outline of the proof of Theorem 1$'$}.
\medskip
In Proposition 2.1, we will show that for $b>\sqrt{\pi/2}$, there
exists a unique function $\phi$ satisfying (1.15)--(1.17) and having
smallest period $b$. Moreover, the initial condition $\phi_0$
increases with  $b$. Next write
$$
\beta\ =\ \sqrt{\pi}\,b,\qquad\qquad\qquad f_0\ =\ e^{\phi_0}\ -\ \phi_0,\qquad\qquad\qquad
M\ =\ \frac1{2b}\int_0^b  \phi (1+e^\phi)\,dy.\tag 1.20
$$
Let us emphasize that although we are now using $4$ variables,
$b,\beta,\phi_0,f_0$, each one is an increasing function of any of
the others.  The non-trivial relationship between them is the
differential equation which relates $b$  to $\phi_0$. We are trying
to prove inequality 1.18, which we write as
$$
M\ +\  \frac{\beta^2}3\ -\ \log (4\beta^2)\ +\ 1\
-\ 4\sum_{n=1}^\infty \log\left| 1-e^{-2n\beta^2}\right|\ <\ 0. \tag 1.21
$$
In Proposition 2.4 we show that
$$
\frac{d (\beta M)}{d\beta }\ =\ 1\ -\ f_0.\tag 1.22
$$
We then investigate how $f_0$ behaves as a function of $\beta$, so
that we can estimate the left hand side of (1.21). Set
$$
\align \ep(\beta)\ &=\ \beta^2\ -\ \log(4\beta^2)\ -\ f_0\\ & =\
\frac{d (\beta M)}{d\beta }\ +\ \beta^2\ -\ \log(4\beta^2)\ -\
1.\tag 1.23
\endalign
$$
We will prove the three key estimates, (1.24)--(1.26).   Set
$\beta_1$ to be the value of $\beta$ corresponding to the initial
value $\phi_0=\log 5$.
$$
\ep(\beta)\ >\ 0,\qquad\qquad \text{ for }\quad
\frac{\pi}{2^{1/2}}<\beta\leq \beta_1, \tag 1.24
$$
$$
 \ep(\beta)\ >\ \frac{0.03}{\beta^2}\qquad\qquad\text{ for }\quad
 \beta_1\leq \beta. \tag 1.25
$$
For some $\gamma>0$, we have
$$
\ep(\beta)\ <\ \frac{\gamma}{\beta^2},\qquad\qquad\text{ for }
\frac{\pi}{2^{1/2}}<\beta.\tag 1.26
$$
 Thus $\ep(\beta)$ is integrable.
 For the proof of (1.24), see  Proposition 2.6--Corollary 2.8.
 For the other two inequalities, see Lemma 2.9 and Proposition 2.10.

 In Corollary 2.5, we obtain a
 simple upper bound on $\beta$ in terms of $\phi_0$ which yields
 $$
 \beta_1\ \leq\ 3.8,\qquad\qquad\text{ for }\phi_0\leq\log 5.
 $$
 Hence integrating (1.24), (1.25) from
$\beta$ to infinity yields
$$
\frac1\beta\int_\beta^\infty \ep(\tilde \beta)\,d\tilde \beta\ >\
\frac{0.01}{\beta^2},\qquad\qquad \text{ for
}\quad\beta>\frac\pi{2^{1/2}}.\tag 1.27
$$
Now integrating (1.23) gives
$$
 M\ +\  \frac{\beta^2}3\ -\ \log (4\beta^2)\ +\ 1\ =\ \frac{C}\beta\
 -\
 \frac1\beta\int_\beta^\infty \ep(\tilde \beta)\,d\tilde \beta, \tag 1.28
$$
where $C$ is the constant of integration. In Proposition 2.11 we
rework some of the asymptotic formulas required in the proof of
(1.25)-(1.26) to show that $C=0$. Hence combining this with (1.27)
gives
$$
 M\ +\  \frac{\beta^2}3\ -\ \log (4\beta^2)\ +\ 1\ \leq\
 -\frac{0.01}{\beta^2},\qquad\qquad \text{ for
}\quad\beta>\frac\pi{2^{1/2}}.\tag 1.29
$$
Finally, one can check with a simple numerical calculation that
$$
- 4\sum_{n=1}^\infty \log\left| 1-e^{-2n\beta^2}\right|\ <\
\frac{0.002}{\beta^2}, \tag 1.30
$$
holds at the value $\beta=\pi/2^{1/2}$. But then in Lemma 2.12 we
see that (1.30) must hold at all values $\beta>\pi/2^{1/2}$. Adding
(1.29) and (1.30) gives (1.21), thus completing the proof of Theorem
1$'$.
\medskip

Now we fill in the results stated in the  outline to complete the
proof.
\medskip\medskip\medskip

\noindent{\bf Section 2. Auxiliary Results and
Proofs}.\medskip\medskip

\proclaim{Proposition 2.1} There exists a smooth function $
\psi:\left(\sqrt{ \pi/2}\ ,\ \infty\right)\times\Bbb R\to \Bbb R $
such that for each fixed $b\in (\sqrt{\pi/2},\infty)$ the function
$$
\phi(y)=\psi(b,y)
$$
satisfies (1.15)--(1.17),  has smallest period $b$, and attains its
minimum value at $y=b/2$. Moreover, writing
$$
f(\phi)\ =\ e^\phi\ -\ \phi,
$$
 $\phi$ is also characterized by having period $b$ and satisfying the following two
conditions:
 $$
 \phi(-y)=\phi(y),\tag 2.1
 $$
$$
y\ =\   \frac1{4\sqrt{\pi}} \int_{\phi(y)}^{\phi_0}
\frac{d\phi}{\sqrt{f_0-f(\phi)}},\qquad\qquad  y\in (0,b/2). \tag
2.2
$$
Furthermore, the map
$$
b\ \mapsto\ \phi_0\ =\ \psi(b,0)
$$
is smooth from the interval $(\sqrt{\pi/2},\infty)$ onto the
interval $(0,\infty)$, and
$$
\frac{db}{d\phi_0}\ >\ 0.
$$
\endproclaim

\remark{Remarks} 1. Every  solution of (1.15)--(1.17) has the form
$$
\phi(y)\ =\ \psi(b/n,y),
$$
for some $n\in \Bbb N$.
\medskip
\noindent 2. By making the change of variables $h=e^\phi$, and
$d\al=e^\phi\,dy$,  we can transform  equation (2.4) to
$$
\frac{d^2 h}{d\al^2}\ =\ 8\pi\left( 1-\frac1h \right). \tag 2.3
$$
Now $d\al$ is a measure of the change in area, and in some respects
it turns out to be more natural to analyze (2.3) than (1.15).
However, we will require a delicate estimate on the relationship
between $b$ and $\phi_0$, and although we work with the variable $h$
at some points, there are places where it is better to work with
(1.15). (For example Proposition 2.4.)
\endremark
\medskip

\noindent{\bf Proof of Proposition 2.1}.  This result is standard
and is part of the standard theory of  ordinary differential
equations, see for example [A] and [Chi]. We give the proof here to
set up notation for later. For $\phi\in\Bbb R$, set
$$
 f(\phi)=e^{\phi}-\phi.
$$
We start by constructing the inverse of $f$.  Indeed, $f$ maps $\Bbb
R$ onto $[1,\infty)$, and for each $f_1\in[1,\infty)$ there exist at
most two solutions of the equation $f(\phi)=f_1$, given by
 $\phi=\phi_*(f_1)$ and $\phi=\phi^*(f_1)$, where
$$
\phi_*(f_1)\leq 0,\qquad\qquad \phi^*(f_1)\geq 0.\tag 2.4
$$

For $\phi_0>0$, we consider the initial value problem
$$
\align & \frac{d^2\phi}{dy^2}\ =\ 8\pi\left( 1-e^\phi\right), \tag 2.5\\
& \phi(0)\ =\ \phi_0.\tag 2.6\\
& \frac{d\phi}{dy}(0)\ =\ 0, \tag 2.7\\
\endalign
$$

Set
$$
f_0=f(\phi_0).\tag 2.8
$$
 Multiplying (2.5) by $d\phi/dy$ and integrating from $y=0$ gives
$$
\left(\frac{d\phi}{dy}\right)^2\ =\ 16\pi\left(e^{\phi_0}-\phi_0\ -\
(e^\phi-\phi)\right)\ =\ 16\pi (f_0-f(\phi)).\tag 2.9
$$
 Hence
$$
\frac{dy}{d\phi}\ =\ \frac{\pm 1}{4\sqrt{\pi}} \
\frac{1}{\sqrt{f_0-f(\phi)}}.\tag 2.10
$$
Set
$$
\ell\ =\ \ell(\phi_0)\ :=\ \frac1{4\sqrt{\pi}}
\int_{\phi_*(f_0)}^{\phi_0} \frac{d\phi}{\sqrt{f_0-f(\phi)}}.\tag
2.11
$$
Then the function $\phi(y)$, assuming it exists, satisfies
$$
y=I(\phi(y))\qquad\text{ for }0\leq y\leq \ell,\qquad\text{ where }\qquad I(z)\ =\   \frac1{4\sqrt{\pi}} \int_{z}^{\phi_0}
\frac{d\phi}{\sqrt{f_0-f(\phi)}}.\tag
2.12
$$
Defining $\phi$ to be the inverse of the function $I$, we find that
$\phi$ is decreasing and smooth on $(0,\ell)$ and it  extends to be
continuously differentiable  on $[0,\ell]$, and satisfies
$$
\phi(0)=\phi_0,\qquad \phi(\ell)=\phi_*(f(\phi_0)),\qquad\qquad
\frac{d\phi}{dy}(0)=\frac{d\phi}{dy}(\ell)=0.
$$
We now extend $\phi$  to $[-\ell, \ell]$  by requiring that it is
even, that is $\phi(-y)=\phi(y)$, and then we extend it to $\Bbb R$
by requiring that it is periodic with period $2\ell$. The result is
an even, continuously differentiable, periodic function on $\Bbb R$
whose smallest period is $2\ell$, and which is smooth on $\Bbb
R\setminus 2\ell\Bbb Z$ and satisfies (2.5) there, and which attains
its maximum value at $y=0$ and its minimum value at $y=\ell$. Now by
the general theorem on the uniqueness and smoothness of solutions to
ordinary differential equations, this solution $\phi$ is smooth and
satisfies (2.5) everywhere on $\Bbb R$. Moreover by the smooth
dependence of solutions to ordinary differential equations on the
initial conditions, we see that defining
$$
\eta(\phi_0,y)\ =\ \phi(y),\qquad\qquad\text{ where }\phi\text{
satisfies (2.5)--(2.7)},
$$
 then  $\eta\in C^\infty((0,\infty)\times\Bbb R)$.
 The final step  is to show that the function
 $$
 \phi_0\ \to\ b=2\ell(\phi_0)
 $$
 is smooth and bijective from  $(0,\infty)$ to
 $(\sqrt{\pi/2},\infty)$, with
 $$
 \frac{db}{d\phi_0}>0,
 $$
 so the inverse function
 $$
 b\ \to\ \phi_0(b)
 $$
 is smooth and bijective from  $(\sqrt{\pi/2},\infty)$
 to $(0,\infty)$.  We then define the function $\psi$  by
 $$
 \psi(b,y)\ =\ \eta(\phi_0(b),y).
 $$
  Proposition 2.1 is thus reduced to the following.

\proclaim{Proposition 2.2} The function $\beta:[1,\infty)\to
[0,\infty)$ defined by
$$
\beta(f_0)\ :=\ \frac12 \int_{\phi_*(f_0)}^{\phi^*(f_0)}
\frac{d\phi}{\sqrt{f_0-f(\phi)}}\tag 2.13
$$
is a smooth function mapping $(1,\infty)$ bijectively onto
$\left(\pi/\sqrt{2}\ ,\ \infty\right)$, with
$$
\frac{d\beta}{df_0}\ >\ 0,\qquad\qquad \text{ on }\qquad(1,\infty).
$$
\endproclaim

\noindent{\bf Proof}.  See [Chi] for a general proof of this result.
See also [ChiJ]. We include the proof here to develop  properties of
the variable $J=j^*+j_*$ which will be useful later on. To reduce
the need for notation, it is convenient to work with physical
variables rather than functions. (To be more precise, we suppose
that there is a fixed underlying ``physical" space which we don't
need to specify. A variable is then a continuous function defined on
this space.) We suppose then that $\phi$ is a variable taking values
in $\Bbb R$, and $f$ and $h$ are variables related to $\phi$ by
$$
f=e^\phi-\phi,\qquad\qquad h=e^\phi,\qquad\qquad \phi=\log
h,\qquad\qquad f=h-\log h. \tag 2.14
$$
The variables $f$ and $h$ take values in $[1,\infty)$ and
$(0,\infty)$ respectively.   Given a value for $f$, we write
$\phi^*\geq 0$ and $\phi_*\leq 0$  for the two corresponding values
for $\phi$ and set
$$
h^*=e^{\phi^*},\qquad\qquad h_*=e^{\phi_*}.\tag 2.15
$$
 When $f=1$ we have $\phi^*=\phi_*=0$ and $h^*=h_*=1$.  For other
values of $f$ the values of $\phi^*$ and $\phi_*$ are distinct. Then
making a change of variables,
$$
\align \frac12 \int_{\phi_*(f_0)}^{\phi^*(f_0)}
\frac{d\phi}{(f_0-f)^{1/2}} \ &=\ \frac12\int_1^{f_0}
\frac1{(f_0-f)^{1/2}} \left( \frac1{e^{\phi^*(f)}-1}\ +\
\frac1{1-e^{\phi_*(f)}}\right)\,df\\
&=\ \frac12\int_1^{f_0} \frac1{(f_0-f)^{1/2}} \left(
\frac1{h^*(f)-1}\ +\ \frac1{1-h_*(f)}\right)\,df. \tag 2.16
\endalign
$$
We will now analyze the Jacobian factor in (2.16) and   modify it to
obtain a positive monotonically increasing function of $f$.

\proclaim{Lemma 2.3} Define variables $j^*$ and $j_*$  by
$$
j^*\ =\ \frac1{h^*-1}\ -\ \frac1{(2(f-1))^{1/2}},\qquad \qquad j_*\
=\ \frac1{1-h_*}\ -\ \frac1{(2(f-1))^{1/2}}.\tag 2.17
$$
Then

(a) As $f\to 1$,
$$ j^*\ \to\
-\frac13,\qquad \qquad j_*\ \to\ \frac13.
$$

(b)  The variables $j^*$ and $j_*$ are increasing with $f$, indeed
 $$
 \frac{dj^*}{df}>0,\qquad\frac{dj_*}{df}>0,\qquad \qquad \text{ for }f>1,
 $$
 and
 $$
\frac{dj^*}{df}\ =\ O((f-1)^{-1/2}),\qquad\frac{dj_*}{df}\ =\ O((f-1)^{-1/2}),
\qquad\qquad\text{ as }f\to 1.
$$

(c) As functions of the variable $f$, the variables  $j^*$ and $j_*$
are concave.  More precisely,

$$
\frac{d^2j^*}{df^2}<0,\qquad\frac{d^2
j_*}{df^2}<0,\qquad\qquad\text{ for }f>1.
$$

(d) The variable
$$
j^*+j_*\ =\ \frac1{h^*-1}\ +\ \frac1{1-h_*}\ -\
\left(\frac2{f-1}\right)^{1/2}\tag 2.18
$$
satisfies
$$
\frac{d(j^*+j_*)}{df}\ >\ 0,\qquad\qquad \frac{d^2(j^*+j_*)}{df^2}\
<\ 0,\qquad\qquad \text{ for }f>1,
$$
and
$$
 (j^*+j_*)\to 0,\qquad\qquad \frac{d(j^*+j_*)}{df}\ =\ O((f-1)^{-1/2}),\qquad\qquad\text{ as }f\to 1.
$$
(e)
$$
0\ <\ j^*+j_*\ <\ 1,\qquad\qquad \text{ when }\ \ f>1.
$$
\endproclaim
\medskip

\noindent{\bf Proof of Lemma 2.3}. Clearly (d) follows from (a), (b)
and (c).  Moreover, see from (2.18) that $j^*+j_*\to 1$,  as $f\to
\infty$,  so (e) follows from (d).
\medskip

\noindent(a) Dealing with the variables $j^*$ and $j_*$
simultaneously, note that as $h\to 1$, we have
$$
\multline
 \frac1{|h-1|}\ -\ \frac1{(2(f-1))^{1/2}}\ =\
 \frac1{|h-1|}\ -\ \frac1{(2(h-1-\log(1-(1- h))))^{1/2}}\\ =\
\frac1{|h-1|}\ -\ \frac1{\left( (1-h)^2\ +\ 2(1-h)^3/3\ +\
2(1-h)^4/4\ +\ \dots\right)^{1/2}}
\endmultline
\tag 2.19
$$
$$ \to\ \cases -1/3 \qquad \text{ as }h\downarrow 1,\\ \ \ 1/3\qquad \text{
as }h\uparrow 1.\endcases
$$
(b)  We need to show that
$$
\frac{d}{df}\left( \frac1{|h-1|}\ -\ \frac1{(2(f-1))^{1/2}}\right)\
>\ 0,
\qquad \qquad \text{ when }\ \ f>1.\tag 2.20
$$

Note that
$$
\frac{dh}{df}\ =\ \frac{h}{h-1}.\tag 2.21
$$
 We thus compute the
sign of the derivative
$$
\multline \frac{d}{df}\left(\frac1{|h-1|}\ -\
\frac1{(2(f-1))^{1/2}}\right) \ =\
\frac{-\sign(h-1)}{(h-1)^2}\frac{dh}{df}
\ +\ \frac{1}{(2(f-1))^{3/2}}\\
=\ \frac{-h}{|h-1|^3} \ +\ \frac{1}{(2(f-1))^{3/2}}.
\endmultline
$$
Hence (2.20) will follow if we can show that
$$
\frac{|h-1|^3}h\ > \ (2(f-1))^{3/2},\qquad\qquad \text{ for }h\neq
1,
$$
equivalently
$$
h^{-2/3}(h-1)^2 \ >\ 2(f-1),\qquad\qquad \text{ for }h\neq 1.\tag
2.22
$$
But this indeed holds, since
$$
h^{-2/3}(h-1)^2 \ -\ 2(f-1)\tag 2.23
$$
equals zero at $h=1$, and
$$
\frac{d}{df} \left( h^{-2/3}(h-1)^2 \ -\ 2(f-1)\right)\ =\  \frac{
2h^{-2/3}(2h+1)}3\ -\ 2\ >\ 0\qquad \text{ for }f>1.
$$
Indeed,
$$
h^{-2/3}(2h+1)\ >\ 3 \qquad\qquad\text{ for }f>1,
$$
as one can easily check by cubing both sides or differentiating once
more with respect to $f$.  The behavior of the derivative as $f\to
1$ is obtained with a Taylor expansion as in (2.19).
\medskip

\noindent (c) We compute
$$
\multline \frac{d^2}{df^2}\left(\frac1{|h-1|}\ -\
\frac1{(2(f-1))^{1/2}}\right) \ =\ \frac{d}{df}\left(
\frac{-h}{|h-1|^3} \ +\ \frac{1}{(2(f-1))^{3/2}}\right)\\
=\  \frac{h(2h+1)}{|h-1|^5}\ -\ \frac{3}{(2(f-1))^{5/2}}.
\endmultline
$$
In order to show that this is negative, we need to show
$$
\frac{(2(f-1))^{5/2}}3\ <\ \frac{|h-1|^5}{h(2h+1)}\qquad\qquad\text{
for }f>1,
$$
or equivalently we need to show
$$
2(f-1)\ <\ 3^{2/5}(h(2h+1))^{-2/5}(h-1)^2\qquad\qquad\text{ for
}f>1.\tag 2.24
$$
Now defining
$$
\tau\ =\ 3^{2/5}(h(2h+1))^{-2/5}(h-1)^2\ -\ 2(f-1),\tag 2.25
$$
we see that $\tau$ vanishes at $f=1$.  Differentiating with respect
to $f$ we get
$$
\frac{d\tau}{df}\ =\  \frac{ 2 \cdot 3^{2/5}}{5}(2h^2+h)^{-7/5}
h\,(6h^2+8h+1)\ -\ 2,\tag 2.26
$$
which also vanishes at $f=1$.  To show that this is positive, we
compute
$$
\frac{5}{2\cdot 3^{2/5}}\frac{d^2\tau}{df^2}\ =\
\frac{2h^2(6h^2-2h+1)}{5 (2h^2+h)^{12/5}}\ >\ 0\qquad\qquad\text{ for
}h>0.\qquad\qquad\qed
$$
\medskip\medskip

Now  we can complete the proof of Proposition 2.2. Introduce the
function $J:[1,\infty)\to \Bbb R$ such that
$$
 j^*+j_*=J(f).
$$
We see that $\beta$ is smooth by fixing $c$ with $1<c<f_0$ and
writing
$$
\multline
 \beta(f_0)\ =\ \frac12\int_1^{f_0} \frac1{(f_0-f)^{1/2}}
\left(\frac2{(2(f-1))^{1/2}}\ +\ J(f)\right)\,df\ =\
\frac{\pi}{2^{1/2}}\ +\ \frac12\int_1^{f_0}
\frac{J(f)}{(f_0-f)^{1/2}} \,df\\
 =\ \frac{\pi}{2^{1/2}}\ +\ \frac12\int_1^{c}
\frac{J(f)}{(f_0-f)^{1/2}} \,df\ +\ \frac12\int_0^{f_0-c}
\frac{J(f_0-f)}{f^{1/2}} \,df.
\endmultline
\tag 2.27
$$
Since $J$ is smooth away from  $1$, both integrals on the right can
be differentiated repeatedly in $f_0$, and we see $\beta$ is smooth
in $f_0$. Differentiating and letting $c\to 0$ gives
$$
\frac{d\beta(f_0)}{df_0}\ =\  \frac12\int_0^{f_0-1}
\frac{J'(f_0-f)}{f^{1/2}} \,df\ >\ 0.\qquad\qquad\qquad\qed
$$

\medskip\medskip

\noindent Our mission is to compute the quantity $M$ in terms of
$\beta$, and we will prove (1.22) relating $M$ to $f_0$. We rescale
the function $\phi$ to have period $2$, by taking the solution
$\psi$ from Proposition 2.1, and setting
$$
\rho(b,s)\ =\ \psi(b,  bs/2),
$$
so that for $b$ fixed, the function $s\mapsto\rho(b,s)$ is even, and
attains its maximum value at $s=0$, and
$$
 \frac{\partial^2\rho}{\partial s^2}\ =\ 2
\beta^2(1-e^{\rho}),\tag 2.28
$$
The solution $\rho$ is  a smooth function of $(\beta,s)$, and we are
interested in the quantity $M$, defined in (1.20).  Setting
$f_0=f(\phi_0)=e^{\phi_0}-\phi_0$, we have from the definition
(1.20), the symmetry of $\phi$, and (2.10),
$$
M\ =\ \frac1{ b}\int_0^{b/2} \phi(1+e^\phi)\,dy\ =\ \frac12\int_0^1
\rho(1+e^\rho)\,ds\ =\ \frac1{4\beta}\int_{\phi_*(f_0))}^{\phi_0}
\frac{\phi(1+e^{\phi})}{(f_0-f(\phi))^{1/2}}\,d\phi.\tag 2.29
$$
\proclaim{Proposition 2.4} (a)
$$
\frac{dM}{d\beta}\ =\ \frac1{\beta} \int_0^1 \rho(1-e^\rho)\,ds.
$$
(b)
$$
 \frac{d(\beta M)}{d\beta}\ =\ \frac12\int_0^1 \rho(3-e^\rho)\,ds\ =\ 1-f_0.
$$
\endproclaim

\noindent{\bf Proof}. (a) We differentiate (2.28) to obtain
$$
\frac{\partial^2}{\partial s^2}\frac{\partial \rho}{\partial \beta}\ =\
4\beta(1-e^{\rho})\ -\ 2\beta^2 \frac{\partial \rho}{\partial \beta}
e^\rho.\tag 2.30
$$
Integrating (2.30) we get
$$
\int_0^1 \frac{\partial \rho}{\partial \beta} e^\rho\,ds\ =\
0.\tag 2.31
$$
Hence
$$
\frac{dM}{d\beta}\ =\   \frac12\int_0^1
\frac{d\rho}{d\beta}(1+e^\rho+\rho e^\rho)\,ds\ =\ \frac12\int_0^1
\frac{d\rho}{d\beta}(1-e^\rho+\rho e^\rho)\,ds.\tag 2.32
$$
However, integrating (2.30) against $\rho$, we get
$$
 \int_0^1 \frac{\partial
\rho}{\partial \beta}\frac{\partial^2\rho}{ds^2}\,ds\ =\
4\beta\int_0^1 \rho(1-e^{\rho})\,ds\ -\ 2\beta^2\int_0^1
\frac{\partial \rho}{\partial \beta} \rho e^\rho\,ds.
$$
Hence using the equation (2.28), we get
$$
 2\beta^2\int_0^1 \frac{\partial\rho}{\partial \beta}(1-e^\rho) \,ds\ =\
 - 2\beta^2\int_0^1
\frac{\partial \rho}{\partial \beta} \rho e^\rho\,ds\ +\
4\beta\int_0^1 \rho(1-e^{\rho})\,ds.
$$
Hence
$$
\frac12 \int_0^1 \frac{\partial\rho}{\partial \beta}(1-e^\rho+\rho
e^\rho) \,ds\ =\ \frac1{\beta}\int_0^1 \rho(1-e^{\rho})\,ds.
$$
Combining this with (2.32) gives (a).
\medskip\medskip

\noindent(b) The first equality follows directly from (a). For the
second, we multiply (2.28) by $d\rho/ds$ and integrating as in
(2.9), to get
$$
\left(\frac{\partial\rho}{\partial s}\right)^2\ =\
4\beta^2(f_0+\rho-e^\rho).
$$
But then
$$
\multline \frac12\int_0^1 \rho(1-e^\rho)\,ds\ =\ \frac1{4\beta^2}
\int_0^1 \rho \frac{\partial^2 \rho}{\partial s^2}\,ds\ =\
-\frac1{4\beta^2}\int_0^1\left(\frac{\partial \rho}{\partial
s}\right)^2\,ds\\ =\ -\int_0^1 (f_0+\rho-e^\rho)\,ds\ =\ 1-f_0-
\int_0^1 \rho\,ds.\qquad\qed
\endmultline
$$
\medskip
\medskip

 \proclaim{Corollary 2.5}
$$
\beta\ \leq\ \frac{\pi}{2^{1/2}}\ +\ (f_0-1)^{1/2}.
$$
\endproclaim

 \noindent{\bf Proof}. From (2.27) and Lemma 2.3 (e), we have
 $$
 \beta(f_0)\ =\ \frac{\pi}{2^{1/2}}\ +\ \frac12\int_1^{f_0}
\frac{J(f)}{(f_0-f)^{1/2}} \,df\ \leq\ \frac{\pi}{2^{1/2}}\ +\
\frac12\int_1^{f_0} \frac{1}{(f_0-f)^{1/2}} \,df\ =\
\frac{\pi}{2^{1/2}}\ +\ (f_0-1)^{1/2}.
 $$
\medskip\medskip

 Our task now is to work towards the estimate in (1.24).  This
 inequality can be checked quite carefully using mathematica, but
 we give a  concise analytic proof with minimal computation.

\proclaim{Proposition 2.7} Given a constant $\la>0$, define
functions $V,W:[1,\infty)\to \Bbb R$ by
$$
\align
V(f)\ &=\ \frac{\pi}{2^{1/2}}\ +\ \la(f-1)^{3/2},\\
 W(f)\ &=\ V(f)^2\
-\ \log(4 V(f)^2)\ -\ f.
\endalign
$$
Suppose that for $f_1>1$ fixed,  there exists  $\la$ such that
$$
\align
  0 <\la<& \frac{2J(f_1)}{3(f_1-1)}, \tag a\\
   W(f_1) &>0,  \tag b \\
 W' (f_1)&<0. \tag c
 \endalign
 $$
Then writing $\beta=\beta(f_0)$ for the function defined in (2.13),
we have
$$
\beta^2-\log(4 \beta^2)-f_0 > 0,\qquad\qquad 1<f_0<f_1.
$$
\endproclaim

\noindent{\bf Proof.}  First we show that $W''(f)>0$ for $f\geq 1$.
Indeed, note that $V(f)>1$ and $V''(f)>0$, and
$$
W'(f)\ =\ 2 \left( V(f)-\frac1{V(f)}\right) V'(f)\ -\ 1,
$$
$$
W''(f)\ =\ 2 \left( 1 +\frac1{V(f)^2}\right)( V' (f))^2\ +\ 2 \left(
V(f) -\frac1{V(f)}\right)V''(f)\ >\ 0.
$$
\medskip

\noindent Next note that $W''>0$  combined with (c) shows that $W$
is decreasing on $[1,f_1]$, and this combined with (b) shows that
$W(f_0)>0$ for $1<f_0<f_1$.
\medskip

\noindent Now we show that for $1<f_0<f_1$ we have
$\beta(f_0)>V(f_0)$. Indeed, comparing the concave function $J(f)$
with the linear function, we get
$$
J(f)\ >\ \frac{J(f_1)}{f_1-1}(f-1),\qquad\qquad 1< f <f_1.
$$
Substituting $t=(f-1)/(f_0-1)$ we get
$$
\multline \beta(f_0)\ =\ \frac{\pi}{2^{1/2}}\ +\ \frac12\int_1^{f_0}
\frac{J(f)}{(f_0-f)^{1/2}} \,df \\  \geq\ \frac{\pi}{2^{1/2}}\ +\
\frac{J(f_1)}{2(f_1-1)}\int_1^{f_0} \frac{f-1}{(f_0-f)^{1/2}} \,df \
=\    \frac{\pi}{2^{1/2}}\ +\
\frac{J(f_1)(f_0-1)^{3/2}}{2(f_1-1)}\int_0^{1} \frac{t}{(1-t)^{1/2}} \,dt \\
=\ \frac{\pi}{2^{1/2}}\ +\ \frac{2J(f_1)(f_0-1)^{3/2}}{3(f_1-1)}\
>\ \frac{\pi}{2^{1/2}}\ +\ \la(f_0-1)^{3/2}\ =\ V(f_0).
\endmultline
$$
Hence we have
$$
\beta^2\ -\ \log(4\beta^2)\ -\ f_0\ >\ V(f_0)^2\ -\ \log(4V(f_0))\
-\ f_0\ =\ W(f_0)\ >\ 0.
$$
\medskip\medskip

\remark{Remark} It will be useful to know the formula
$$
h_*\ =\ \sum_{j=0}^\infty \frac{(j+1)^{j-1}}{j!}\ e^{-(j+1)f},
$$
although we will not prove it or depend on it. \endremark\medskip

\proclaim{Lemma 2.8} For $\la=0.098$, the conditions of Proposition
2.7 are satisfied for $f_1=5-\log 5$.
\endproclaim
\medskip

\noindent{\bf Proof of Lemma 2.8}.  Step 1: For $h=5$, write
$5_*=h_*$. Then numerical calculation shows that
$$
5-\log(5)= f(5)\ <\ f(0.034)=0.034-\log(0.034).
$$
Hence
$$
5_*>0.034,
$$
and
$$
\multline (j+j_*)(f(5))\ =\ \frac1{5-1}\ +\ \frac1{1-5_*}\ -\
\left(\frac2{5-\log 5-1}\right)^{1/2} \ \geq\ \frac1{4}\ +\
\frac{1}{1-0.034}\ -\ \left(\frac2{4-\log 5}\right)^{1/2}\\ =\
0.3705...,
\endmultline
$$
and the right hand term in Proposition 2.6 (a) is
$$
\frac{2(j+j_*)(f(5))}{3(4-\log 5)}\ =\ 0.1033....\ >\ 0.098.
$$
\medskip

\noindent Step 2:
$$
V(f(5))\ =\ \frac{\pi}{2^{1/2}}\ +\ \la(4-\log 5)^{3/2}\ =\
2.583664...
$$
so
$$
2.58366\ <\ V(f(5))\ <\ 2.584.
$$
Hence at the value $f_1=5-\log 5$ we have
$$
W\ =\ V^2\ -\ \log(4V^2)\ -\ (5-\log 5)\ \geq\
2.58366^2-\log(4\times 2.58366^2)-(5-\log 5)\ =\ 0.00002...\ >\ 0,
$$
while
$$
W'\ =\ 3\la \left( V -\frac1V\right)(4-\log 5)^{1/2}\ -\ 1\ <\
0.294\left( 2.584-\frac1{2.584}\right)(4-\log 5)^{1/2}\ -\ 1\, =\,
-0.001...\, <\, 0. \qed
$$
\medskip\medskip

\proclaim{Corollary 2.8} Set $\beta_1=\beta(5-\log 5)$. For
$\beta=\beta(f_0)$, we set
$$
\ep(\beta)\ =\ \beta^2\ -\ \log(4\beta^2)\ -\ f_0.
$$
Then
$$
\ep(\beta)\ >\ 0,\qquad\qquad \text{ if }\quad
\frac{\pi}{2^{1/2}}<\beta< \beta_1.
$$
\endproclaim

\noindent{\bf Proof}.  By Lemma 2.7, if $1<f_0<f_1=5-\log 5$, then
$\ep(\beta)>0$.

\medskip
\medskip
\medskip

\noindent Now we will investigate more precisely how $b(f_0)$
depends on $f_0$ as $f_0\to\infty$.  We will use the fact that we
only have positive Taylor coefficients in the expansion
$$
(1-x)^{-1/2}\ =\ \sum_{0}^\infty \gamma_k x^k,  \qquad\qquad
\gamma_k=\frac{(2k)!}{2^{2k} (k!)^2}\sim \frac1{\sqrt{\pi k}}.
$$
From (2.5) and the fact that $\phi$ is even and periodic with period
$b$, we have
$$
\int_0^{b/2} (1-e^\phi)\,dy\ =\ 0.
$$
Hence using (2.10) and setting $h=e^\phi$, we get
$$
\beta\ =\ \sqrt{\pi}b\ =\ 2\sqrt{\pi}\int_0^{b/2}e^\phi\,dy \ =\
\frac12 \int_{\phi_*(f_0)}^{\phi^*(f_0)} \frac{e^\phi
\,d\phi}{\sqrt{f_0-f(\phi)}}\ =\ \frac12 \int_{h_*(f_0)}^{h^*(f_0)}
\frac{dh}{\sqrt{f_0-(h-\log h)}}.\tag 2.33
$$
Using the notation of (2.14), (2.15) and writing $h_0=h^*(f_0)$ and
$h_{0*}=h_*(f_0)$, so $f_0=h_0-\log h_0=h_{0*}-\log h_{0*}$, and
setting $t=1-h/h_0$, we get

$$
\align \beta\ &=\ \frac1{2} \int_{h_{0*}}^{h_0} \left(h_0-h+\log
\frac{h}{h_0}\right)^{-1/2}\,dh
 \\  &=\ \frac{h_0}{2} \int_{0}^{1-h_{0*}/h_0} \left(h_0 t\ +\ \log
(1-t)\right)^{-1/2}\,dt\\
&=\ \frac{h_0}{2(h_0-1)^{1/2}} \int_{0}^{1-h_{0*}/h_0}t^{-1/2}
\left(1-\frac{-\log (1-t)-t}{(h_0-1)t}\right)^{-1/2}\,dt \tag 2.34\\
&=\ \frac{h_0}{2(h_0-1)^{1/2}}
\sum_{k=0}^\infty\frac{\gamma_k}{(h_0-1)^k}\int_{0}^{1-h_{0*}/h_0}
t^{-1/2}\left(\frac{-\log (1-t)-t}{t}\right)^{k}\,dt\\
&=\
\frac{h_0}{(h_0-1)^{1/2}}\sum_{k=0}^\infty\frac{\gamma_k\,\mu_k(1-h_{0*}/h_0)}{(h_0-1)^k},
\tag 2.35
\endalign
$$
where the series converges by monotone convergence, and
$$
\mu_k(\tau)\ =\ \frac12\int_0^\tau t^{-1/2}\left(\frac{-\log
(1-t)-t}{t}\right)^{k}\,dt.\tag 2.36
$$
Clearly $\mu_k(\tau)$ is an increasing function of  $\tau$ which is
strictly positive for $\tau\in(0,1]$. Now
$$
\beta^2\ =\ \frac{h_0^2}{h_0-1}\sum_{k=0}^\infty
\frac{\nu_k(1-h_{0*}/h_0)}{(h_0-1)^k},\qquad\qquad\qquad
\nu_k(\tau)\ =\ \sum_{j=0}^k \gamma_j\, \gamma_{k-j}\,
\mu_j(\tau)\,\mu_{k-j}(\tau).\tag 2.37
$$
Clearly $\nu_j(\tau)$ is also positive.  It is easy to compute
$$
\mu_0(\tau)\ =\ \tau^{1/2},\qquad\qquad \nu_0(\tau)\ =\ \tau.
$$
 Hence just taking the first term in
(2.37) gives
$$
\beta^2\ >\ \frac{h_0^2 \nu_0(1-h_{0*}/h_0)}{h_0-1}\ =\
\frac{h_0(h_0-h_{0*})}{h_0-1}\ >\ h_0.\tag 2.38
$$
Now applying the Mean Value Theorem to the function $w\mapsto w-\log
w$, we have
$$
\align \beta^2-\log(\beta^2)\ -\ (h_0-\log h_0)\ &\geq\
\frac{h_0-1}{h_0}(\beta^2-h_0)\ =\ \sum_{k=0}^\infty
\frac{h_0\ \nu_k(1-h_{0*}/h_0)}{(h_0-1)^k}\ -\ h_0\ +\ 1\\
&\geq\ 1\ -\  h_{0*}\ +\ \frac{h_0\ \nu_1(1-h_{0*}/h_0)}{h_0-1}\ +\
\frac{h_0\ \nu_2(1-h_{0*}/h_0)}{(h_0-1)^2}.\tag 2.39
\endalign
$$
\proclaim{Lemma 2.9} If $h_0\geq 5$ then
$$
1\ -\  h_{0*}\ +\ \frac{h_0\ \nu_1(1-h_{0*}/h_0)}{h_0-1}\ >\ 2\log
2,\tag a
$$
and
$$
\frac{h_0\ \nu_2(1-h_{0*}/h_0)}{(h_0-1)^2}\ >\
\frac{0.03}{\beta^2}\tag b
$$
so
$$
\ep(\beta)\ =\ \beta^2-\log(4\beta^2)\ -\ f_0\ >\
\frac{0.03}{\beta^2}\tag c
$$
\endproclaim

\noindent{\bf Proof}.  (a) Now evaluating (2.36) for $k=1$,
$$
\mu_1(\tau)\ =\ 2\log(1+\tau^{1/2})\ -\ \tau^{1/2}\ +\
(\tau^{-1/2}-1)\log(1-\tau),\tag 2.40
$$
so  we have
$$
 \nu_1(\tau)\ =\ \mu_0(\tau)\mu_1(\tau)\ =\ 2\tau^{1/2}\log(1+\tau^{1/2})\ -\ \tau\
+\ (1-\tau^{1/2})\log(1-\tau).\tag 2.41
$$
We will estimate the terms on the right hand side. Since
$\log(1-\tau)<0$, we have
$$
(1-\tau^{1/2})\log(1-\tau) \ >\ (1-\tau)\log(1-\tau),
$$
and by the convexity of the logarithm we have $\log(1+x)> x\log 2$
for $0<x<1$,  and so
$$
\tau^{1/2}\log(1+\tau^{1/2})\ >\ \tau\log 2.
$$
Hence substituting these inequalities into (2.41),
$$
\nu_1(\tau)\ >\ \tau(2\log 2-1)\ +\ (1-\tau)\log(1-\tau).
$$
 and writing
$1-\tau=h_{0*}/h_0$  we have
$$
\multline 1\ -\  h_{0*}\ +\ \frac{h_0\ \nu_1(1-h_{0*}/h_0)}{h_0-1}\
\\>\ 1\ -\ h_{0*}\ +\ \frac{h_0}{h_0-1}\left((2\log 2\ -\ 1)\left(1\ -\
\frac{h_{0*}}{h_0}\right)\
+\ \frac{h_{0*}}{h_0}\log\frac{h_{0*}}{h_0}\right)\\
=\  2\log 2\ +\ \frac{(2\log 2-1)\ +\ 2h_{0*}(1-\log 2)\ -\
h_{0*}(h_0+\log(h_0/h_{0*}))}{h_0-1},
\endmultline
$$
and so (a) holds provided
$$
2\log 2-1\ +\ 2h_{0*}(1-\log 2)\ -\
h_{0*}(h_0+\log(h_0/h_{0*}))\ \geq\ 0,\qquad\qquad\text{ for }h_0>5,
\tag 2.42
$$
which certainly follows if we can show
$$
h_{0*}(h_0+\log h_0+\log1/h_{0*})\ <\ (2\log
2-1),\qquad\qquad\text{ for } h_0>5.\tag 2.43
$$
We first remark that
$$
0.035-\log 0.035\ <\ 5-\log 5,
$$
and hence if $h_0>5$, then
$$
h_{0*}\ <\ 0.035\ <\ \exp(-1).
$$
But then for $h_{0*}<0.035$, we have that $h_{0*}\log(1/h_{0*})$
increases with $h_{0*}$ and hence decreases with $h_0$.  Moreover,
$$
h_{0}-1\ >\ 4\ >\ 1-h_{0*},
$$
so the functions
$$
h_0\ \mapsto\  h_{0*}h_0
$$
and
$$
h_0\ \mapsto\ h_{0*}\log h_0
$$
are also decreasing with $h_0$, as can be checked by differentiating
with respect to $f_0$.  For example
$$
\frac{d(h_{0*} h_0)}{df_0}\ =\ h_{0*} h_0\left(
\frac1{h_{0}-1}-\frac1{1-h_{0*}}\right)\ <\ 0.
$$
Hence the left hand side of (2.43) is decreasing with $h_0$, and so
bounded above by
$$
0.035(5+\log 5+\log 1/0.035)\ =\ 0.34866...\ <\ 0.386294..\ =\ 2\log
2-1,
$$
and (2.43) holds, so (a) holds.
\medskip\medskip

\noindent (b) Now for $\tau>0$,
$$
\nu_2(\tau)\ =\ \frac{3 \mu_0(\tau)\mu_2(\tau)}4\ +\
\frac{(\mu_1(\tau))^2}4\ >\ \frac{(\mu_1(\tau))^2}4.
$$
Hence for $h_0\geq 5$, we have
$$
h_{0*}<0.035
$$
and
$$
\frac{(\mu_1(1-h_{0*}/h_0))^2}4\ \geq\ \frac{(\mu_1(1-0.035/5))^2}4\
=\ 0.0340...\ >\ 0.03.
$$
Hence
$$
\frac{h_0\ \nu_2(1-h_{0*}/h_0)}{(h_0-1)^2}\ \geq\
\frac{0.03}{h_0-1}\ >\ \frac{0.03}{\beta^2}.
$$

\noindent (c) Follows by substituting (a) and  (b) into (2.39).\qed
\medskip
\medskip

\proclaim{Proposition 2.10}
$$
\ep(\beta)\ =\ O(\beta^{-2}),\qquad\qquad\text{ as
}\qquad\beta\to\infty.
$$
\endproclaim

\noindent{\bf Proof}. We will prove this by bounding the error when
we approximate the series in (2.35) by the partial sums. Indeed, we
show that there exists a constant $C(K)$ independent of $h_0$ such
that
$$
\left|\beta\ -\ \frac{h_0}{(h_0-1)^{1/2}}\sum_{k=0}^{K}
\frac{\gamma_k\,\mu_k(1)}{(h_0-1)^k}\right|\ {\leq}\ \frac{C(K)
}{h_0^{K+1/2}},\qquad\qquad\text{ for }h_0>2.\tag 2.44
$$
In fact, what we show is
$$
\left|\beta\ -\ \frac{h_0}{(h_0-1)^{1/2}}\sum_{k=0}^{K}
\frac{\gamma_k\,\mu_k(1)}{(h_0-1)^k}\right|\ \leq\ \frac{C(K)(\log
h_0)^{K} }{h_0^{K-1/2}},\qquad\qquad\text{ for }h_0>2. \tag 2.45
$$
By applying (2.45) with $K$ replaced by $K+2$, we get (2.44).
\medskip

\noindent{\bf Notation}. Suppose $h=(h_1,\dots,h_p)$ and
$k=(k_1,\dots,k_q)$ are variables taking values in $U\subset\Bbb
R^p$ and $V\subset R^q$ respectively, and suppose that $F_1$ and
$F_2$ are two functions of $(h,k)$.  Then  we write
$$
F_1\ \leqk\ F_2,
$$
if for every  $k\in V$, there exists a constant $C(k)<\infty$, such
that
$$
F_1(h,k)\ \leq\ C(k) \ F_2(h,k)\qquad\text{ for all }h\in U.
$$
\medskip

Now we prove (2.45). We first remark that
$$
\left|(1-x)^{-1/2}\ -\ \sum_{k=0}^{K-1} \gamma_k x^k\right|\ \leqK\
(1-x)^{-1/2}\,x^{K},\qquad\qquad \text{ for }0\leq x<1.
$$
Hence from (2.34), writing $t=1-h/h_0$, we have that for $h_0>2$,
$$
\multline \left|\beta\ -\ \frac{h_0}{(h_0-1)^{1/2}}\sum_{k=0}^{K-1}
\frac{\gamma_k\,\mu_k(1-h_{0*}/h_0)}{(h_0-1)^k}\right| \\ \leqK\
\frac{h_0}{(h_0-1)^{K+1/2}}\int_{0}^{1-h_{0*}/h_0}
t^{-1/2}\left(1-\frac{-\log
(1-t)-t}{(h_0-1)t}\right)^{-1/2}\left(\frac{-\log
(1-t)-t}{t}\right)^{K}\,dt\\
=\ \frac{1}{(h_0-1)^{K}} \int_{h_{0*}}^{h_0} \left(h_0-\log
h_0-(h-\log h)\right)^{-1/2}\left(\frac{-\log
(1-t)-t}{t}\right)^{K}\,dh.
\endmultline
\tag 2.46
$$
We split into two cases. The function
$$
\frac{-\log (1-t)-t}{t}\ =\
\sum_{k=1}^\infty\frac{t^k}{k+1},\qquad\qquad 0< t<1,
$$
is increasing with $t$, so decreasing with $h$.  Hence  for $h_0>h$
and  $h_0>2$, we have
$$
\frac{-\log (1-t)-t}{t}\ =\ \frac{\log(h_0/h)}{1-h/h_0}\ -\ 1 \
<\ \cases \frac83 \log h_0 \qquad & h>1/h_0,\\
  -\frac83\log h & h\leq 1/h_0.\endcases
$$
Hence the right hand side of (2.46) is bounded up to a constant
$C(K)$ by
$$
\align   \frac{(\log h_0)^{K}}{(h_0-1)^{K}} \int_{h_{0*}}^{h_0}
&\left(h_0-\log
h_0-(h-\log h)\right)^{-1/2}\,dh\\
&+\ \frac{1}{(h_0-1)^{K}} \int_{h_{0*}}^{1/h_0} \left(h_0-\log
h_0-(h-\log h)\right)^{-1/2}(-\log h)^{K}\,dh.\tag 2.47
\endalign
$$
Using Corollary 2.5,  for  $h_0>2$, the first term in (2.47) is
equal to
$$
\frac{2(\log h_0)^K\beta}{(h_0-1)^{K}}\ \leqK\ \frac{(\log
h_0)^K}{(h_0-1)^{K}}\left(\frac\pi{2^{1/2}}+(f_0-1)^{1/2}\right)\
\leqK \ \frac{ (\log h_0)^K}{(h_0-1)^{K-1/2}}.
$$
To bound the second term in (2.47), we  change variables to
$f=h-\log h$ to get the bound
$$
  \frac{1
}{(h_0-1)^{K}} \int_{\log h_0+1/h_0}^{f_0}
\left(f_0-f\right)^{-1/2}(-\log h_*(f))^{K}\frac{h_*(f)}{1-h_*(f)}\,df.\tag 2.48
$$
But
$$
-\log h_*(f)\ =\ f\ -\ h_*(f)\ <\ f\ +\ 1,
$$
and for $f>\log(2/\log 2)$ we have
$$
h_*(f)\ \leq\ 2e^{-f}.
$$
Hence  (2.48) is bounded up to a constant $C(K)$ by
$$
\frac{1}{(h_0-1)^{K}} \int_0^{f_0} \left(f_0-f\right)^{-1/2}e^{-f}
f^K\,df.
$$
But the integral here is uniformly bounded in $f_0$, so the second
term in (2.47) is bounded up to $C(K)$ by
$$
\frac{1 }{(h_0-1)^{K}}.
$$
So far we have bounded the left hand side of (2.46) by the right
hand side of (2.45). To complete the proof of (2.45) we just have to
show that for $h_0>2$,
$$
\left|\mu_k(1-h_{0*}/h_0)-\mu_k(1)\right|\ \leqkK\ \frac{1}{h_0^K}.
$$
However, the left hand side equals
$$
\multline \left|\,\frac12\int_{1-h_{0*}/h_0}^1
t^{-1/2}\left(\frac{-\log (1-t)-t}{t}\right)^{k}\,dt\,\right|\
\leqk\ \int_0^{h_{0*}/h_0} |\log s|^k\,ds\\ \leqk\ \frac{h_{0*}
}{h_0}|\log h_0-\log h_{0*}|^k \ =\  \frac{h_{0*} }{h_0}( h_0-
h_{0*})^k\ \leqk\  e^{-f_0} h_0^{k-1}\ \leqkK\ \frac{1}{h_0^K}.
\endmultline
$$
This completes the proof of (2.45). From this we get from this the
asymptotic formula
$$
\beta^2\qquad \sim\qquad \frac{h_0^2}{h_0-1}\sum_{k=0}^\infty
\frac{\nu_k(1)}{(h_0-1)^k},
$$
where $\nu_k$ is defined in (2.37), in the sense that for $h_0>2$,
$$
\left|\beta^2\ -\ \frac{h_0^2}{h_0-1}\sum_{k=0}^{K}
\frac{\nu_k(1)}{(h_0-1)^k}\right|\ \leqK\ \frac1{h_0^K}.
$$
Thus
$$
\beta^2\ =\ \frac{h_0^2}{h_0-1}\left( 1\ +\ \frac{2\log 2-1}{h_0-1}\
\right)\ +\ O(h_0^{-1})\ =\ h_0\ +\ 2\log 2\ +\ O(h_0^{-1}).\tag
2.49
$$
From this we see that
$$
 \beta^2-\log(\beta^2)\ -\
(h_0-\log h_0)\ -\ 2\log 2\ =\  O(h_0^{-1})\ =\ O(\beta^{-2}).
$$
This completes the proof of Proposition
2.10.\qquad\qquad\qquad\qquad\qquad\qquad\qquad\qquad\hfill \qed
\medskip\medskip\medskip

\proclaim{Proposition 2.11}
$$
M\ =\  -\frac{\beta^2}3\ +\ \log(4\beta^2)\ -\ 1\ +\ O(\beta^{-2})
\qquad\text{ as }\qquad \beta\to \infty.
$$
\endproclaim
\medskip

\noindent{\bf Proof}.  From (2.29), we have
$$
\multline M\ =\ \frac1{4\beta} \int_{\phi_*(f_0)}^{\phi^*(f_0)}
\frac{\phi (e^{\phi}+1)}{\sqrt{f_0-f(\phi)}}\,d\phi\\ =\
\frac1{3\beta} \int_{\phi_*(f_0)}^{\phi^*(f_0)} \frac{(\phi-\log
h_0) e^{\phi}}{\sqrt{f_0-f(\phi)}}\,d\phi\ + \frac{\log h_0}{3\beta}
\int_{\phi_*(f_0)}^{\phi^*(f_0)} \frac{
e^{\phi}}{\sqrt{f_0-f(\phi)}}\,d\phi\ +\ \frac1{12\beta}
\int_{\phi_*(f_0)}^{\phi^*(f_0)} \frac{\phi
(3-e^{\phi})}{\sqrt{f_0-f(\phi)}}\,d\phi\\ =\ \frac1{3\beta}
\int_{\phi_*(f_0)}^{\phi^*(f_0)} \frac{(\phi-\log h_0)
e^{\phi}}{\sqrt{f_0-f(\phi)}}\,d\phi \ +\ \frac{2\log h_0}3\ +\
\frac{1-f_0}3.
\endmultline
\tag 2.50
$$
The third line here follows from (2.33) and the second equality in
Proposition 2.4(b). Now we change variables to $h=e^\phi$ so
$f=e^\phi-\phi=h-\log h$, and set $h_0=h^*(f_0)$ and
$h_{0*}=h_*(f_0)$. Then define
$$
N\ :=\ \frac1{2} \int_{\phi_*(f_0)}^{\phi^*(f_0)} \frac{(\phi-\log h_0)
e^{\phi}}{\sqrt{f_0-f(\phi)}}\,d\phi\ =\ \frac1{2}
\int_{h_{0*}}^{h_0} \frac{\log h-\log h_0}{\sqrt{f_0-f}}\,dh.\tag
2.51
$$
 We follow the argument of (2.34)-(2.35) with $\beta$ replaced by (2.51) to
get
$$
\multline N\ =\ \frac{h_0}{2(h_0-1)^{1/2}}
\int_{0}^{1-h_{0*}/h_0}(\log (1-t))\,t^{-1/2}\left(1-\frac{-\log
(1-t)-t}{(h_0-1)t}\right)^{-1/2}\,dt\\
=\ \frac{h_0}{(h_0-1)^{1/2}}
\sum_{k=0}^\infty\frac{\gamma_k\,\kappa_k(1-h_{0*}/h_0)}{(h_0-1)^k},
\endmultline
$$
where
$$
\kappa_k(\tau)\ =\ \frac12\int_0^\tau (\log(1-t))\,
t^{-1/2}\left(\frac{-\log (1-t)-t}{t}\right)^{k}\,dt.
$$
Moreover, following the proof of (2.44)-(2.45), we conclude that for
$h_0>2$,
$$
\left|N\ -\ \frac{h_0}{(h_0-1)^{1/2}}
\sum_{k=0}^{K-1}\frac{\gamma_k\,\kappa_k(1)}{(h_0-1)^k}\right| \
\leqK\ \frac1{h_0^{K-1/2}}.
$$
Now
$$
\kappa_0(1)=2\log 2-2,
$$
and so in particular, using (2.49),
$$
N\ =\
\kappa_0{h_0^{1/2}} \ +\ O( h_0^{-1/2})\ =\ (2\log 2-2)\beta\ +\
O(\beta^{-1}),\qquad\text{ as } \beta\to \infty.
$$
Substituting this into (2.50) and using (2.49), we see that as
$\beta\to\infty$ we have
$$
\multline M\  =\ \frac{ 2(\log 4-2)}{3} \ +\frac{2\log (\beta^2)}3\
+\ \frac{-\beta^2+ \log (4\beta^2)
+1}3\ +\ O(\beta^{-2})\\
=\  -\frac{\beta^2}3\ +\ \log(4\beta^2)\ -\ 1\ +\ O(\beta^{-2}).
\endmultline
$$
\medskip\medskip

\proclaim{Lemma 2.12} Suppose that $C>0$ and $\beta_1>1/\sqrt{2}$
are constants and that the formula
$$
 - 4\sum_{n=1}^\infty
\log\left( 1-e^{-2n\beta^2}\right)\ <\ \frac{C}{\beta^2},\tag 2.52
$$
holds for $\beta=\beta_1$.  Then it holds for all
$\beta\geq\beta_1$.
\endproclaim

\noindent{\bf Proof}. Define
$$
\omega(\beta)\ =\ - 4\sum_{n=1}^\infty \log\left(
1-e^{-2n\beta^2}\right),\qquad\qquad \beta>0,
$$
and
$$
\psi(\beta)\ =\ \frac{C}{\beta^2}\ -\ \omega(\beta).
$$
Then $\omega$ is positive and smooth, and
$$
-\omega'(\beta)\ =\ 16 \beta \sum_{n=1}^\infty
\frac{n}{1-e^{-2n\beta^2}}\ >\ 16 \beta \sum_{n=1}^\infty
\frac{1}{1-e^{-2n\beta^2}}\ >\ 4\beta \omega(\beta).
$$
 Suppose that (2.52) fails, that is $\psi(\beta)\leq 0$,  for some
$\beta_2>\beta_1$.  Then we can choose $\beta_2>\beta_1$ minimal
such that this is the case, and clearly $\psi(\beta_2)=0$.  But then
$$
\psi'(\beta_2)\ =\ \frac{2C}{\beta_2^3}\ -\ \omega'(\beta_2)\ =\
\frac{2\omega(\beta_2)}{\beta_2}\ -\ \omega'(\beta_2)\ \geq\
\omega(\beta_2)\left( \frac{-2}{\beta_2}\ +\ 4\beta_2\right).
$$
But $\beta_2>1/\sqrt{2}$, so the right hand side is positive and so
$\psi(\beta)<0$  for some $\beta$ with $\beta_1<\beta<\beta_2$,
which is a contradiction.\qquad\qquad\qed

\medskip\medskip\medskip
 \noindent{\bf Appendix. Explicit formulas for the flat torus and
the round sphere}.
\medskip

\proclaim{Lemma A.1} Let $T=\Bbb C/\Lambda$ be a torus of area $1$,
where $\Lambda$ is a lattice, and let $u$ and $v$ be the generators
of the dual lattice $\Lambda^*$ and set $z=v/u$.  Then for the flat
metric $g_0$ on $T$,
$$
\trace \Delta_{g_0}^{-1}\ =\ -
\frac{\log 2\pi}{2\pi} \ -\ \frac{\log(|\eta(z)|^4/|u|^2)}{4\pi},
\tag A.1
$$
where the Dedekind eta function $\eta$ is defined by
$$
\eta(z)\ =\ e^{\pi iz/12}\prod_{n=1}^\infty (1-e^{2\pi i nz}).\tag
A.2
$$
On the other hand,
$$
\trace\Delta^{-1}_{S^2,1}\ =\ -\frac{\log \pi}{4\pi}\ -\
\frac1{4\pi},\tag A.3
$$
and so
$$
\trace\Delta_{g_0}^{-1}\ -\ \trace\Delta_{S^2,1}^{-1}\ =\
\frac{1}{4\pi}\left(-\log(|\eta(z)|^4/|u|^2)\ -\ \log 4\pi\ +\
1\right).\tag A.4
$$
When $\Lambda$ has generators $(1/b,a+bi)$ with $a,b\in \Bbb R$, we
can choose $(u,v)=(-i/b, b-ai )$ and then (A.4) becomes (1.12).
\endproclaim

\remark{Remark}.  The quantity $\log (|\eta(z)|^4/|u|^2)$ was shown
in [OsPS] to be maximized at the hexagonal torus, for which
$$
\log (|\eta(z)|^4/|u|^2)\ =\ -1.0335...
$$
Hence the hexagonal torus minimizes $\trace\Delta^{-1}$ among flat
tori of a given area.
\endremark

\noindent{\bf Proof}.  Now
$$
\Lambda^*\ =\ \{\mu\in\Bbb C:\Re(\bar \mu \lambda)\in\Bbb Z\text{ for all }\lambda\in \Lambda\}.
$$
The eigenfunctions of the Laplacian on $ T=\Bbb C/\Lambda$ have the
form
$$
f(z)\ =\ e^{2\pi i \Re(\bar\mu \lambda)}  ,\qquad\qquad \text{ for }\mu\in \Lambda^*.
$$
The corresponding eigenvalue is  $(2\pi)^2 |\mu|^2$. Consider the
Epstein zeta function
$$
Z_T(s)\ =\ \sum_{\mu\in \Lambda ^*-0} \frac1{(2\pi
|\mu|)^{2s}}.
$$
{\it Kronecker's First Limit Formula} states that
$$
(2\pi)^{2s} Z_T(s)\ =\ \frac{\pi}{s-1}\ +\ 2\pi\left( -\Gamma'(1)\ -\ \log 2\ -\ \log|\eta(z)|^2\right)
\ +\ O(s-1).
$$
Hence
$$
Z_T(s)\ =\  \frac1{4\pi(s-1)}\ +\ Z_T^1\ +\ O(s-1),\qquad\qquad
Z_T^1\ =\ \frac1{2\pi}\left( -\Gamma'(1)\ -\ \log (4\pi)\ -\ \log|\eta(z)|^2\right).
$$
But $Z_T^1$ is a different regularization of the trace of
$\Delta^{-1}$, and it can be shown that this differs from our Green
function regularization $\trace\Delta^{-1}$ by a universal constant:
$$
\trace\Delta^{-1}_{g_0}\ =\ Z_T^1\ +\ \frac{\log
2}{2\pi}\ +\ \frac{\Gamma'(1)}{2\pi} . \tag A.5
$$
see  [M2], [S1], [S2], or [O2] (A.6).  Evaluating (A.5) we get
(A.1).

Formula (A.3) is well known.  Indeed, on the round $2$-sphere of
area $4\pi$ given by $x^2+y^2+z^2=1$, the Green function $G(p,q)$
can be written in terms of the distance $r$ from $p$ to $q$, as
$$
G(p,q)\ =\ -\frac1{2\pi} \log |\sin r/2|\ -\ \frac1{4\pi}.
$$
This gives the Robin mass
$$
m_{S^2,4\pi}\ =\ \frac{\log 2}{2\pi}\ -\ \frac1{4\pi},
$$
and combining this with  (1.3) gives
$$
m_{S^2,1}\ =\ m_{S^2,4\pi}\ -\ \frac{\log 4\pi}{4\pi}\  =\ -\frac{\log \pi}{4\pi}\ -\
\frac1{4\pi}.
$$
\medskip
The author is extremely grateful to the referee for pointing out
several  results related to this work and providing helpful
comments.
\medskip

\centerline{References}
\medskip

\roster

\item"[A]" V. Arnold: Ordinary differential equations.  Universitext. Springer-Verlag, Berlin,
2006.

\item"[B]" W. Beckner: Sharp Sobolev inequalities on the sphere and the Moser-Trudinger inequality.
{\it Annals of Math.} {\bf 138} (1993), 213-242.

\item"[CL]" E. Carlen and M. Loss: Competing symmetries, the logarithmic HLS inequality
and Onofri's inequality on $S^n$.  {\it Geometric and Functional
Analysis} {\bf 2} (1992) 90--104.

\item"[Ch]" S.-Y. A. Chang: Conformal invariants and partial differential equations.
{\it Bull. Amer. Math. Soc.} {\bf  42} (2005),  365--393.

\item"[Chi]" C. Chicone: The monotonicity of the period function for
planar Hamiltonian vector fields.  {\it J. Differential Equations}
{\bf 69} (1987), no. 3, 310--321.

\item"[ChiJ]"  C. Chicone and M. Jacobs: Bifurcation of limit cycles from
quadratic isochrones.  {\it J. Differential Equations}  {\bf  91}
(1991), no. 2, 268--326.

\item"[Chiu]" P. Chiu: Height of flat tori. {\it Proc. Aner. Math.
Soc.,} {\bf 125} 723--730, (1997)

\item"[DS1]"  P. Doyle and J. Steiner: Spectral invariants and playing hide and seek on surfaces.
{\it Preprint}.

\item"[DS2]"  P. Doyle and J. Steiner: Blowing bubbles on the torus.
{\it Preprint}.

\item"[LL1]" C.-S. Lin and M. Lucia:
Uniqueness of solutions for a mean field equation on torus. {\it J.
Differential Equations}  {\bf 229} (2006),  no. 1, 172--185.

\item"[LL2]" C.-S. Lin and M. Lucia: One-dimensional symmetry of
periodic minimizers for a mean field equation. {\it Ann. Sc. Norm.
Super. Pisa Cl. Sci.} (5)  {\bf 6}  (2007),  no. 2, 269--290.

\item"[M1]"  C. Morpurgo:  The logarithmic Hardy-Littlewood-Sobolev inequality
and extremals of zeta functions on $S^n$. {\it Geom. Funct. Anal.}
{\bf 6} (1996),  146--171.

\item"[M2]"  C. Morpurgo: Zeta functions on $S\sp 2$. Extremal Riemann
surfaces (San Francisco, 1995), 213--225, {\it Contemp. Math.}, {\bf
201}, Amer. Math. Soc., Providence, RI, (1997).

\item"[M3]" C. Morpurgo:
 Sharp inequalities for functional integrals and traces of conformally invariant operators.
 {\it Duke Math. J.} {\bf 114} (2002),  477--553.

\item"[O1]" K. Okikiolu:   Hessians of spectral zeta functions. {\it Duke Math. J.}
{\bf  124} (2004), 517--570.

\item"[O2]" K. Okikiolu:  Extremals for Logarithmic HLS inequalities on compact
manifolds. {\it GAFA} (To Appear.)

\item"[OW]" K. Okikiolu and C. Wang:  Hessian of the zeta function for the Laplacian on forms.
{\it Forum Math.} {\bf  17} (2005),  105--131.

\item"[On]" Onofri: E. On the positivity of the effective action in a theory of random surfaces.
{\it Comm. Math. Phys.} {\bf  86} (1982),  321--326.

\item"[OsPS]"  B. Osgood, R. Phillips and P. Sarnak: Extremals of determinants of Laplacians.
{\it J. Funct. Anal.} {\bf 80} (1988),  148--211.

\item"[S1]"  J. Steiner: {\it Green's Functions, Spectral Invariants, and a Positive Mass on Spheres}.
Ph. D. Dissertation, University of California San Diego, June 2003.

\item"[S2]" J. Steiner:  A geometrical mass and its extremal properties for metrics on $S\sp 2$.
{\it Duke Math. J.}  {\bf 129} (2005), 63--86.

\item"[T]" A. Terras, {\it Harmonic analysis and symmetric spaces and applications I.}
 Springer-Verlag, Berlin, 1988.

\endroster
\medskip\medskip

\centerline{Kate Okikiolu}

\centerline{University of California, San Diego}

\centerline{okikiolu\@math.ucsd.edu}

\end
\end